\documentclass[a4paper,12pt]{amsart}
\usepackage{amssymb,amscd,amsmath,a4wide,graphicx,stmaryrd,fullpage,setspace,microtype,xr,ulem,textcomp,csquotes}
\SetSymbolFont{stmry}{bold}{U}{stmry}{m}{n}
\usepackage[pagebackref=true]{hyperref}
\usepackage[all]{xy}
\usepackage[T1]{fontenc}

\title{Enriched diagrams of topological spaces over locally contractible enriched categories}

\author{Philippe Gaucher}

\address{Universit\'e de Paris, IRIF, CNRS, F-75013 Paris, France}

\urladdr{http://www.irif.fr/{\~{}}gaucher} 



\swapnumbers

\newcommand{\C}{\mathcal{C}}
\newcommand{\D}{\mathcal{D}}
\newcommand{\K}{\mathcal{K}}

\newcommand{\W}{\mathcal{W}}
\newcommand{\F}{\mathcal{F}}

\newcommand{\p}\times

\newtheorem*{thmN}{Theorem}

\newtheorem{thm}{Theorem}[section]
\newtheorem{prop}[thm]{Proposition}

\newtheorem{cor}[thm]{Corollary}

\newtheorem*{conv}{Convention}
\newtheorem{defn}[thm]{Definition}
\newtheorem{nota}[thm]{Notation}

\newcommand{\bd}{\begin{defn}}
	\newcommand{\ed}{\end{defn}}
\newcommand{\bp}{\begin{prop}}
	\newcommand{\ep}{\end{prop}}
\newcommand{\bth}{\begin{thm}}
	\renewcommand{\eth}{\end{thm}}
\newcommand{\bpf}{\begin{proof}}
	\newcommand{\epf}{\end{proof}}
\newcommand{\bc}{\begin{cor}}
	\newcommand{\ec}{\end{cor}}

\newcommand{\fL}[1]{\ar@{->}[ll]_-{#1}}
\newcommand{\fR}[1]{\ar@{->}[rr]^-{#1}}
\newcommand{\fRr}[1]{\ar@{->}[rrr]^-{#1}}
\newcommand{\fD}[1]{\ar@{->}[dd]_-{#1}}
\newcommand{\fU}[1]{\ar@{->}[uu]^-{#1}}
\newcommand{\f}[2]{\ar@{->}[#1]|{#2}}
\newcommand{\ff}[2]{\ar@2{->}[#1]|{#2}}
\newcommand{\frr}[1]{\ar@{->}[rrrr]^-{#1}}

\newcommand{\fl}[1]{\ar@{->}[l]_-{#1}}
\newcommand{\fr}[1]{\ar@{->}[r]^-{#1}}
\newcommand{\fd}[1]{\ar@{->}[d]_-{#1}}
\newcommand{\fu}[1]{\ar@{->}[u]^-{#1}}

\renewcommand{\top}{{\mathbf{Top}}}

\newcommand{\iso}{\cong}

\newcommand{\ot}{\otimes}

\renewcommand{\leq}{\leqslant}
\renewcommand{\geq}{\geqslant}

\newcommand{\dgr}[1]{\top^{#1_0}}
\newcommand{\dgrQ}[1]{\top_Q^{#1_0}}

\def\cartesien{%
	\ar@{-}[]+R+<6pt,-2pt>;[]+RD+<6pt,-6pt>%
	\ar@{-}[]+D+<2pt,-6pt>;[]+RD+<6pt,-6pt>%
}
\def\cocartesien{%
	\ar@{-}[]+L+<-6pt,+2pt>;[]+LU+<-6pt,+6pt>%
	\ar@{-}[]+U+<-2pt,+6pt>;[]+LU+<-6pt,+6pt>%
}
\def\hocartesien{%
	\ar@{-}[]+R+<6pt,-2pt>;[]+RD+<6pt,-6pt>_{h}%
	\ar@{-}[]+D+<2pt,-6pt>;[]+RD+<6pt,-6pt>%
}
\def\hococartesien{%
	\ar@{-}[]+L+<-6pt,+2pt>;[]+LU+<-6pt,+6pt>_{h}%
	\ar@{-}[]+U+<-2pt,+6pt>;[]+LU+<-6pt,+6pt>%
}

\newcommand{\brm}[1]{\rm{\mathbf{#1}}}

\newcommand{\set}{{\brm{Set}}}

\newcommand{\ttop}{{\brm{TOP}}}

\DeclareMathOperator{\id}{Id}

\DeclareMathOperator{\Obj}{Obj}
\DeclareMathOperator{\Mor}{Mor}

\newcommand{\liminj}{\varinjlim}
\newcommand{\limproj}{\varprojlim}
\newcommand{\dcat}{{\mathbf{Cat}}}

\newcommand{\thin}[1]{\mathrm{Thin}(#1)}

\makeatletter
\def\varholim@#1#2{%
	\vtop{\m@th\ialign{##\cr
			\hfil$#1\operator@font holim$\hfil\cr
			\noalign{\nointerlineskip\kern1.5\ex@}#2\cr
			\noalign{\nointerlineskip\kern-\ex@}\cr}}%
}
\def\holimproj{%
	\mathop{\mathpalette\varholim@{\leftarrowfill@\textstyle}}\nmlimits@
}
\def\holiminj{%
	\mathop{\mathpalette\varholim@{\rightarrowfill@\textstyle}}\nmlimits@
}
\makeatother

\DeclareMathOperator{\cell}{{\brm{cell}}}
\DeclareMathOperator{\cof}{{\brm{cof}}}
\DeclareMathOperator{\inj}{{\brm{inj}}}
\newcommand{\ddownarrow}{{\downarrow}}

\newcommand{\adj}[6]{\xymatrix@C=#5em{{#1}\ar@/^#6pt/[r]^{#2} \ar@{}[r]|{\perp} & \ar@/^#6pt/[l]^{#3} {#4}}}

\begin{document}

\begin{abstract} It is proved that the projective model structure of the category of topologically enriched diagrams of topological spaces over a locally contractible topologically enriched small category is Quillen equivalent to the standard Quillen model structure of topological spaces. We give a geometric interpretation of this fact in directed homotopy. 
\end{abstract}

\maketitle

\tableofcontents

\section{Introduction}

\subsection*{Presentation}
In directed homotopy, we have to deal with various mathematical objects encoding execution paths and irreversibility of time. In the frameworks of \textit{$d$-spaces} in Grandis' sense \cite{mg} and of \textit{multipointed $d$-spaces} in the sense of \cite{mdtop} (the latter is a variant of the former), it is considered a topological space of states equipped with a distinguished set of continuous maps $dX$. The set of execution paths is in the case of Grandis' $d$-spaces invariant by reparametrization by the monoid $\mathcal{M}$ of nondecreasing continuous maps from $[0,1]$ to itself preserving the extremities. The set of execution paths of a multipointed $d$-space is only invariant by reparametrization by the group $\mathcal{G}$ of nondecreasing homeomorphisms of $[0,1]$ (the reason for the difference in the axiomatization is that we are still unable to replace everywhere $\mathcal{G}$ by $\mathcal{M}$ to prove all what is known about multipointed $d$-spaces, in particular the results of \cite{model2}). For this introduction, $\mathcal{P}$ is either $\mathcal{G}$ or $\mathcal{M}$. 

The starting point of the paper is the following geometric observation. The monoid $\mathcal{P}$ can be viewed as a one-object category such that $\mathcal{P}$ is the unique set of morphisms. Consider a pair $(X,dX)$ such that $X$ is a topological space and such that $dX$ is a set of continuous maps from $[0,1]$ to $X$ invariant by the action of $\mathcal{P}$. We suppose that $dX$ is equipped with its natural topology making the evaluation maps continuous. This data gives rise to a contravariant diagram $\D^{\mathcal{P}}(X,dX)$ of topological spaces over $\mathcal{P}$  with the only vertex $dX$ and taking $\phi\in \mathcal{P}$ to the mapping $\phi^*:\gamma\mapsto \gamma.\phi$. The limit $\limproj \D^{\mathcal{P}}(X,dX)$ is the space of paths of $dX$ invariant by the action of $\mathcal{P}$. It is equal to the subspace of constant paths of $dX$. The colimit $\liminj \D^{\mathcal{P}}(X,dX)$ is nothing else but the quotient of the space of paths $dX$ by the action of $\mathcal{P}$. More interesting is the interpretation of the homotopy colimit $\holiminj \D^{\mathcal{P}}(X,dX)$. The expected behavior is that this homotopy colimit is weakly homotopy equivalent to $dX$ because calculating the colimit up to homotopy should prevent the identifications up to reparametrization from being made. It is the case if $\mathcal{P} = \mathcal{M}$ by Theorem~\ref{good-behavior}. It turns out that it is not the case for the homotopy colimit $\holiminj \D^{\mathcal{G}}(X,dX)$. For example, if $dX$ is a singleton (i.e. a constant path), then $\holiminj \D^{\mathcal{G}}(X,dX)$ has the homotopy type of $B\mathcal{G}$ by \cite[Proposition~14.1.6 and Proposition~18.1.6]{ref_model2} which is not contractible, although both $\mathcal{G}$ and $dX$ (which is supposed to be here a singleton) are contractible. The main result of this paper is that one possible way to overcome this problem is to work in the enriched setting. The main theorem of this paper is stated now: 

\begin{thmN} (Theorem~\ref{eq-topdgr-top})
Let $\mathcal{P}$ be a topologically enriched small category. Suppose that $\mathcal{P}$ is locally contractible (i.e all spaces of maps $\mathcal{P}(\ell,\ell')$ are contractible). Let $\top$ be the category of $\Delta$-generated spaces. Then the colimit functor from the category $[\mathcal{P},\top]_0$ of topologically enriched functors  and natural transformations to $\top$ induces a left Quillen equivalence between the projective model structure and the Quillen model structure.
\end{thmN}

Note that the particular case where $\mathcal{P}$ has exactly one map between each pair of objects (i.e. each space $\mathcal{P}(\ell,\ell')$ is a singleton) is trivial. In this case, $[\mathcal{P},\top]_0$ is equivalent to $\top$ as a category indeed.

Using this theorem, our example can now be reinterpreted in the enriched setting. The diagram $\D^{\mathcal{P}}(X,dX)$ belongs to $[\mathcal{P}^{op},\top]_0$ because the mapping $\phi\mapsto \phi^*$ from $\mathcal{P}$ to $\ttop(dX,dX)$, where $\ttop(dX,dX)$ is the space of continuous maps from $dX$ to itself, is continuous. Since the inclusion functor $[\mathcal{P}^{op},\top]_0\subset \top^{\mathcal{P}^{op}}$ into the category of all contravariant functors from $\mathcal{P}$ to $\top$ is colimit-preserving and limit-preserving by Proposition~\ref{colim-lim-enriched}, nothing changes concerning the interpretations of $\limproj \D^{\mathcal{P}}(X,dX)$ and $\liminj \D^{\mathcal{P}}(X,dX)$. On the contrary, the behavior of the homotopy colimit is completely different. There exists a cofibrant replacement $\D^{\mathcal{P}}(X,dX)^{cof}$ of $\D^{\mathcal{P}}(X,dX)$ in $[\mathcal{P},\top]_0$ together with a pointwise weak homotopy equivalence $\D^{\mathcal{P}}(X,dX)^{cof} \to \Delta_{\mathcal{P}^{op}} (dX)$ (the constant diagram $\Delta_{\mathcal{P}^{op}} (dX)$ belongs to $[\mathcal{P}^{op},\top]_0$). Since $dX$ is fibrant, we deduce that the canonical map $\liminj \D^{\mathcal{P}}(X,dX)^{cof} \to dX$ is a weak homotopy equivalence because the colimit functor is a left Quillen equivalence. It means that in the enriched setting, $\holiminj \D^{\mathcal{P}}(X,dX)$ always has the homotopy type of $dX$. Everything behaves as if we were in the trivial case above.

\subsection*{Motivation} Theorem~\ref{eq-topdgr-top} as well as Theorem~\ref{lefproper-ok} are of general interest for all mathematicians studying the algebraic topology of enriched diagrams of $\Delta$-generated spaces. A conjectural generalization of Theorem~\ref{eq-topdgr-top} is even suggested in Section~\ref{conclusion}. 

My own interest for this theorem comes from what follows. I need to force a cofibrant replacement to behave in a very specific way in a model category (the model category of \textit{Moore flows}) which will be introduced in a subsequent paper. Everything boils down to fixing the behavior of the homotopy colimit. There are two ways of fixing this behavior: Replacing $\mathcal{G}$ by $\mathcal{M}$ and using Theorem~\ref{good-behavior} (but I am still unable to replace $\mathcal{G}$ by $\mathcal{M}$ in all proofs of \cite{model2} and I need the results of these papers to reach the goal explained below) or tweaking the notion of Moore flow (there are several types of Moore flows indeed) to use Theorem~\ref{eq-topdgr-top}. The ultimate goal is to prove that there is a zig-zag of Quillen equivalences between the model category of \textit{multipointed $d$-spaces} \cite{mdtop} and the model category of \textit{flows} \cite{model3} \cite{leftdetflow} with one of the model categories of Moore flows in the middle: it is a work in progress. The only known relation between these two model categories is \cite[Theorem~7.5]{mdtop} which can be reformulated as follows. There exists a functor from multipointed $d$-spaces to flows which is neither a left nor a right adjoint and such that the total left derived functor in the sense of \cite{HomotopicalCategory} induces an equivalence of categories between the homotopy categories. I still do not know how to prove the latter theorem by replacing $\mathcal{G}$ by $\mathcal{M}$ in the definition of multipointed $d$-space.

\subsection*{Outline of the paper}
Section~\ref{reminder} collects the notations and some useful facts about locally presentable categories and model categories which are used in this paper. Section~\ref{base} gives a very important example of locally presentable base in the sense of \cite{EnrichedSketch}. Section~\ref{Qmixed} recalls some facts about the category of $\Delta$-generated spaces, the \textit{Quillen model structure}, the \textit{Cole-Str{\o}m model structure} and the so-called \textit{mixed model structure}. It culminates in the proof that the mixed model structure is accessible. Section~\ref{enricheddiag} introduces the material of enriched diagrams of topological spaces. Some elementary facts which are used in the next sections are proved or recalled. Section~\ref{homotopytheory} introduces two model structures, the projective one, which is combinatorial, on the category of enriched diagrams using the Quillen model structure, and the injective one, which is only accessible, on the category of enriched diagrams using the mixed model structure. This section discusses the interactions between the two model structures. The existence of the projective model structure is a straightforward consequence of Moser's work \cite{MoserLyne}. We are also able to prove that this model structure is left proper in Theorem~\ref{lefproper-ok} (it is right proper because all objects are fibrant). The latter result is not a consequence of Moser's work. It is based on the results of a previous work \cite[Appendix~A]{mdtop} which partially relies on a weakening of the notion of closed $T_1$-inclusion introduced by Dugger and Isaksen in \cite[p~686]{hocolimfacile}. Section~\ref{mainthm} proves the main theorem of the paper. Section~\ref{conclusion} adds a comment about the monoid of nondecreasing continuous maps from $[0,1]$ to itself preserving the extremities and another one mentioning Shulman's work \cite{shuhomotopycolim} about enriched homotopical categories. Finally, an appendix proves two particular cases of Theorem~\ref{eq-topdgr-top} using \cite{shuhomotopycolim}. 

\subsection*{Acknowledgments}
I thank Tim Campion for Theorem~\ref{timthm} and Tyler Lawson for Theorem~\ref{bm-contractible}. I also thank Asaf Karagila for \cite{2948666} even if his contribution is not necessary anymore in this new version (it suffices in the proof of Theorem~\ref{timthm} to choose a regular cardinal $\mu >\lambda$ such that $\mu\rhd \lambda$ instead of $\mu^\lambda=\mu$) and Tim Porter for valuable email discussions. I thank the anonymous referee for the detailed report.

\section{Notations, conventions and prerequisites}
\label{reminder}

It is not required to read any paper about directed homotopy to understand this work.  We refer to \cite{TheBook} for locally presentable categories, to \cite{MR2506258} for combinatorial model categories.  We refer to \cite{MR99h:55031} and to \cite{ref_model2} for more general model categories. We refer to \cite{KellyEnriched} and to \cite[Chapter~6]{Borceux2} for enriched categories. \textit{All enriched categories are topologically enriched categories: the word topologically is therefore omitted.} 
\begin{itemize}
\item All categories are locally small except the category $\mathbf{CAT}$ of all locally small categories.
\item $\K$ always denotes a locally presentable category.
\item $\set$ is the category of sets.
\item $\K(X,Y)$ is the set of maps in a category $\K$.
\item $\K_0$ denotes sometimes the underlying category of an enriched category $\K$. Since we will be working in a cartesian closed category of topological spaces, $\K_0$ is nothing else but the category $\K$ with the topology of the space of maps forgotten.
\item $\dcat$ is the category of all small categories and functors between them.
\item $\mathcal{P}$ denotes a nonempty enriched small category.
\item $\K^{op}$ denotes the opposite category of $\K$.
\item $\Obj(\K)$ is the class of objects of $\K$.
\item $\Mor(\K)$ is the category of morphisms of $\K$ with the commutative squares for the morphisms.
\item $\K^I$ is the category of functors and natural transformations from a small category $I$ to $\K$.
\item $[I,\K]$ is the enriched category of enriched functors from an enriched small category to an enriched category $\K$, and $[I,\K]_0$ is the underlying category.
\item $\Delta_I:\K \to \K^I$ is the constant diagram functor.
\item $\varnothing$ is the initial object.
\item $\mathbf{1}$ is the final object.
\item $\id_X$ is the identity of $X$.
\item $g.f$ is the composite of two maps $f:A\to B$ and $g:B\to C$; the composite of two functors is denoted in the same way. 
\item $F\Rightarrow G$ denotes a natural transformation from a functor $F$ to a functor $G$.
\item A subcategory is always isomorphism-closed (replete).
\item $f\boxslash g$ means that $f$ satisfies the \textit{left lifting property} (LLP) with respect to $g$, or equivalently that $g$ satisfies the \textit{right lifting property} (RLP) with respect to $f$.
\item $\inj(\C) = \{g \in \K, \forall f \in \C, f\boxslash g\}$.
\item $\cof(\C)=\{f\mid \forall g\in \inj(\C), f\boxslash g\}$.
\item $\cell(\C)$ is the class of transfinite compositions of pushouts of elements of $\C$.
\item A \textit{cellular} object $X$ of a combinatorial model category is an object such that the canonical map $\varnothing\to X$ belongs to $\cell(I)$ where $I$ is the set of generating cofibrations.
\item A model structure $(\C,\W,\F)$ means that the class of cofibrations is $\C$, that the class of weak equivalences is $\W$ and that the class of fibrations is $\F$ in this order.
\item $(-)^{cof}$ denotes a cofibrant replacement, $(-)^{fib}$ denotes a fibrant replacement.
\item $F\dashv G$ denotes an adjunction where $F$ is the left adjoint and $G$ the right adjoint.
\end{itemize}
We will use the following known facts: 
\begin{itemize}
\item A functor $F:\K\to\mathcal{L}$ between locally presentable categories is a left adjoint if and only if it is colimit-preserving; indeed, any left adjoint is colimit-preserving; conversely, if $F$ is colimit-preserving, $F^{op}$ is limit-preserving; since every locally presentable category is well-copowered by \cite[Theorem~1.58]{TheBook} and has a generator, the opposite category $\K^{op}$ is well-powered and has a cogenerator; Hence the Special Adjoint Functor theorem \cite[Theorem~3.3.4]{Borceux1} states that $F^{op}$ is a right adjoint.
\item A functor $F:\K\to\mathcal{L}$ between locally presentable categories is a right adjoint if and only if it is limit-preserving and accessible by \cite[Theorem~1.66]{TheBook}.
\item A transfinite tower (of length $\lambda$) of $\K$ consists of an ordinal $\lambda$ and a colimit-preserving functor $D$ from $\lambda$ to $\K$; it means that for every limit ordinal $\mu\leq \lambda$, the canonical map $\liminj_{\nu<\mu} D_\nu\to D_\mu$ is an isomorphism. In any model category $\mathcal{M}$, a colimit of a transfinite tower of cofibrations between cofibrant objects is a homotopy colimit. It is due to the fact that the transfinite tower is a diagram over a direct Reedy category and that, in this case, the tower is Reedy cofibrant for the Reedy model structure which coincides with the projective model structure \cite[Theorem~5.2.5]{MR99h:55031}.
\end{itemize}

A weak factorization system $(\mathcal{L},\mathcal{R})$ of a locally presentable category $\K$ is \textit{accessible} if there is a functorial factorization 
\[\xymatrix@1{(A\stackrel{f}\longrightarrow B) \ar@{|->}[r] & (A\stackrel{Lf}\longrightarrow Ef\stackrel{Rf}\longrightarrow B)}\]
with $Lf\in \mathcal{L}$, $Rf\in \mathcal{R}$ such that the functor $E:\Mor(\K)\to \K$ is accessible \cite[Definition~2.4]{GKR18}. Since colimits are calculated pointwise in $\Mor(\K)$, a weak factorization system is accessible if and only if the functors $L:\Mor(\K)\to \Mor(\K)$ and $R:\Mor(\K)\to \Mor(\K)$ are accessible. By \cite[Theorem~4.3]{MR3638359}, a weak factorization system is accessible if and only if it is small in Garner's sense. In particular, every \textit{small} weak factorization system (i.e. of the form $(\cof(I),\inj(I))$ for a set $I$) is accessible. A model structure $(\C,\W,\F)$ on a locally presentable category is \textit{accessible} if the two weak factorization systems $(\C,\W\cap\F)$ and $(\C\cap\W,\F)$ are accessible. Every combinatorial model category is an accessible model category.

We will be using the following characterization of a Quillen equivalence. A Quillen adjunction $F\dashv G:\C\leftrightarrows \D$ is a Quillen equivalence if and only if for all fibrant objects $X$ of $\D$, the natural map $F(G(X)^{cof})\to X$ is a weak equivalence of $\D$ (the functor $F$ is then said to be \textit{homotopically surjective}) and if for all cofibrant objects $Y$ of $\C$, the natural map $Y\to G(F(Y)^{fib})$ is a weak equivalence of $\C$ \cite[Proposition~1.3.13]{MR99h:55031}. If all objects of $\D$ are fibrant, the latter assertion is equivalent to saying that for all cofibrant objects $Y$ of $\C$, the unit of the adjunction $Y\to G(F(Y))$ is a weak equivalence of $\C$.

\section{On locally presentable bases}
\label{base}

\bd \cite[Definition~1.1]{EnrichedSketch} Let $\lambda$ be a regular cardinal. A locally $\lambda$-presenta\-ble base is a symmetric monoidal closed category $\mathcal{V}$ which is locally $\lambda$-presentable and such that 
\begin{itemize}
	\item The unit of the tensor product is $\lambda$-presentable
	\item The tensor product of two $\lambda$-presentable objects of $\mathcal{V}$ is $\lambda$-presentable.
\end{itemize}
A locally presentable base is a locally $\lambda$-presentable base for some regular cardinal $\lambda$. 
\ed

The following theorem is the key fact to establish Corollary~\ref{localbase}. 

%

\bth \cite{TimHelp} (T. Campion) \label{timthm}
Let $\K$ be a locally presentable category. Let $\lambda$ be a regular cardinal. Then there exists a regular cardinal $\mu > \lambda$ such that the binary product of two $\mu$-presentable objects is $\mu$-presentable.  
\eth

\bpf 
Since we have the isomorphisms \[\K(Z,X\p Y)\iso \K(Z,X)\p \K(Z,Y) \iso (\K\p \K)((Z,Z),(X,Y)),\] the functor $(X,Y)\mapsto X\p Y$ is a right adjoint. It is therefore accessible. We choose a large enough regular cardinal $\nu$ such that the functor $(X,Y)\mapsto X\p Y$ is $\nu$-accessible and such that $\nu > \lambda$. We choose a regular cardinal $\mu \rhd \nu$ such that the binary product of two $\nu$-presentable objects is $\mu$-presentable. By \cite[Remark~2.15]{TheBook}, write $A$ as a retract of a $\mu$-small $\nu$-filtered colimit $\overline{A} = \liminj_{k\in K} A_k$ of $\nu$-presentable objects and $B$ as a retract of a $\mu$-small $\nu$-filtered colimit $\overline{B} = \liminj_{\ell\in L} B_{\ell}$ of $\nu$-presentable objects. Let $I=K\p L$ which is $\mu$-small and $\nu$-filtered. The projections $\pi_1:K\p L\to K$ and $\pi_2:K\p L\to L$ are right cofinal. Indeed, for $k\in K$, $k\ddownarrow \pi_1=(k\ddownarrow K)\p L$ is a product of filtered categories, and so filtered itself and therefore nonempty and connected. This implies that $\pi_1$ (and also $\pi_2$) is right cofinal. We obtain the isomorphisms $\overline{A}\iso\liminj_{i\in I} A_{\pi_1(i)}$ and $\overline{B}\iso\liminj_{i\in I} B_{\pi_2(i)}$ by \cite[Theorem~14.2.5(1)]{ref_model2}. We deduce the isomorphism of $\K\p\K$ \[(\overline{A},\overline{B})\iso\liminj_{i\in I} \left(A_{\pi_1(i)},B_{\pi_2(i)}\right).\] Since $I$ is $\nu$-filtered and since the functor $(X,Y)\mapsto X\p Y$ is supposed to be $\nu$-accessible, we obtain the isomorphism 
\[\overline{A}\p \overline{B}\iso  \liminj_{i\in I} \left(A_{\pi_1(i)} \p B_{\pi_2(i)}\right).\] Let $C=\liminj_{j\in J} C_j$ be a $\mu$-filtered colimit. Then we have 
\begin{align*}
&\K(\overline{A}\p \overline{B},C) \\&\iso \limproj_{i\in I} \K\left(A_{\pi_1(i)}\p B_{\pi_2(i)},\liminj_{j\in J} C_j\right)&\hbox{by $\overline{A}\p \overline{B} = \liminj\limits_{i\in I} \left(A_{\pi_1(i)} \p B_{\pi_2(i)}\right)$}\\
&\iso \limproj_{i\in I} \liminj_{j\in J} \K(A_{\pi_1(i)}\p B_{\pi_2(i)},C_j)&\hbox{since $A_{\pi_1(i)}\p B_{\pi_2(i)}$ is $\mu$-presentable}\\
&\iso \liminj_{j\in J} \limproj_{i\in I} \K(A_{\pi_1(i)}\p B_{\pi_2(i)},C_j)&\hbox{since $I$ is $\mu$-small and $J$ is $\mu$-filtered}\\
&\iso \liminj_{j\in J} \K(\overline{A}\p \overline{B},C_j)&\hbox{by $\overline{A}\p \overline{B} = \liminj\limits_{i\in I} \left(A_{\pi_1(i)} \p B_{\pi_2(i)}\right)$.}
\end{align*}
We deduce that $\overline{A}\p \overline{B}$ is $\mu$-presentable, and $A\p B$ as well since it is a retract of $\overline{A}\p \overline{B}$. 
\epf

We deduce the following important example of locally presentable base which is not in \cite{EnrichedSketch}: 

\begin{cor}  (T. Campion) \label{localbase} Let $\K$ be a locally presentable category which is cartesian closed. Then it is a locally presentable base for the closed monoidal structure induced by the binary product.
\end{cor}

\bpf It suffices to start from a regular cardinal $\lambda$ such that $\K$ is locally $\lambda$-presentable and such that the terminal object is $\lambda$-presentable and to apply Theorem~\ref{timthm}. 
\epf

\section{Quillen and mixed model structures of topological spaces}
\label{Qmixed}

The category $\top$ denotes the category of \textit{$\Delta$-generated spaces}, i.e. the colimits of simplices. For a tutorial about these topological spaces, see for example \cite[Section~2]{mdtop}. The category $\top$ is locally presentable (see \cite[Corollary~3.7]{FR}), cartesian closed and it contains all CW-complexes. The internal hom functor is denoted by $\ttop(-,-)$. The forgetful functor from $\top$ to $\set$ is fibre-small and topological. The category $\top$ is a full coreflective subcategory of the category $\mathcal{T\!O\!P}$ of general topological spaces. 

The category $\top$ can be viewed as a category enriched over itself. It is also locally presentable in the enriched sense by \cite[Proposition~2.4]{MoserLyne} (see also \cite[Corollary~7.3]{Kelleyb}). It is tensored and cotensored over itself because $\top$ is cartesian closed: the tensor product is the binary product and the unit is the singleton. A category enriched over $\top$ is called an \textit{enriched category}. As already said in Section~\ref{reminder}, the adjective \enquote{topologically} is omitted because all enrichments in this paper are over $\top$. 

We recall Cole's theorem enabling us to mix model structures.

\bth \cite[Theorem~2.1]{mixed-cole} \label{cole-model-structure0}
Let $(\C_1,\W_1,\F_1)$ and $(\C_2,\W_2,\F_2)$ be two model structures on the same underlying category with $\W_1\subset \W_2$ and with $\F_1\subset \F_2$. Then there exists a unique model structure $(\C_m,\W_m,\F_m)$ such that $\W_m=\W_2$ and $\F_m=\F_1$. Moreover, we have $\C_1\cap\W_1 = \C_m\cap \W_m$ and $\C_2\subset \C_m$.
\eth

\bp \label{cole-model-structure} With the notations of Theorem~\ref{cole-model-structure0}. Suppose that the underlying category $\K$ is locally presentable. Suppose that the weak factorization system $(\C_1\cap\W_1,\F_1)$ is accessible and that the model structure  $(\C_2,\W_2,\F_2)$ is combinatorial. Then the model structure $(\C_m,\W_m,\F_m)$ is accessible.
\ep

\bpf There is the equality of weak factorization systems $(\C_1\cap\W_1,\F_1)=(\C_m\cap\W_m,\F_m)$ by Theorem~\ref{cole-model-structure0}. Thus the right-hand weak factorization system is accessible because the left-hand one is accessible by hypothesis. The other factorization is obtained as follows: first $f$ factors as a composite $f=R_2(f).L_2(f)$ with $L_2(f)\in \C_2$ and $R_2(f)\in \W_2\cap \F_2$. Since $\C_2\subset \C_m$ by Theorem~\ref{cole-model-structure0}, $L_2(f)\in \C_m$. Then $R_2(f)$ factors as a composite $R_2(f)=\ell.k$ with $k\in \C_1\cap\W_1 = \C_m\cap \W_m$ and $\ell\in \F_1=\F_m$. By the 2-out-of-3 property, $\ell\in \W_2=\W_m$. Thus the second factorization is $f=\ell.(k.L_2(f))$. The functor $R_2:\Mor(\K)\to \Mor(\K)$ is accessible since the model structure $(\C_2,\W_2,\F_2)$ is combinatorial by hypothesis. Since $(\C_1\cap\W_1,\F_1)$ is accessible by hypothesis, we deduce that the weak factorization system $(\C_m,\W_m\cap\F_m)$ is accessible.
\epf

The category $\top$ can be equipped with the Quillen model structure $(\C_2,\W_2,\F_2)$ in which the weak equivalences are the weak homotopy equivalences \cite[Section~2.4]{MR99h:55031}. There is another well-known model structure $(\C_1,\W_1,\F_1)$ on $\top$ called the Cole-Str{\o}m model structure. The weak equivalences are the homotopy equivalences; the fibrations are the Hurewicz fibrations; the cofibrations are the strong Hurewicz cofibrations. A general proof of its existence can be found in \cite[Corollary~5.23]{Barthel-Riel}; The monomorphism hypothesis is automatically satisfied because $\top$ is locally presentable \cite[Remark~5.20]{Barthel-Riel}. All topological spaces are fibrant and cofibrant for the Cole-Str{\o}m model structure. By using Proposition~\ref{cole-model-structure}, we obtain the mixed model structure: the weak equivalences are the weak homotopy equivalences and the fibrations are the Hurewicz fibrations. All topological spaces are fibrant for this model structure. The cofibrations (the cofibrant objects resp.) of the mixed model structure are called the \textit{mixed cofibrations} (the \textit{mixed cofibrant objects} resp.). All Quillen cofibrations are mixed cofibrations because $\C_2\subset \C_m$. By \cite[Proposition~3.6]{mixed-cole}, a map $f:A\to X$ is a mixed cofibration if and only if it is a closed Hurewicz cofibration and $f$ factors as a composite $f:A\to X' \to X$ such that the left-hand map is a Quillen cofibration and the right-hand map is a homotopy equivalence. In particular, the cofibrant objects of the mixed model structure are the topological spaces homotopy equivalent to a cofibrant object of the Quillen model structure \cite[Corollary~3.7]{mixed-cole}.

\begin{nota} By convention, $\top_Q$ denotes the category of $\Delta$-generated spaces equipped with the Quillen model structure and $\top_m$ denotes the category of $\Delta$-generated spaces equipped with the mixed model structure. \end{nota}

\begin{conv}
The words cofibration and cofibrant without further precision mean cofibration and cofibrant in $\top_Q$. The words mixed cofibration and mixed cofibrant mean cofibration and cofibrant in $\top_m$. 
\end{conv}

\begin{cor} \label{mixed-accessible}
The model category $\top_m$ is accessible.
\end{cor}

Concerning the Cole-Str{\o}m model structure $(\C_1,\W_1,\F_1)$ of $\top$, it is known that the weak factorization system $(\C_1,\W_1\cap \F_1)$ is not small by \cite[Remark~4.7]{Raptis-Strom}. It is unlikely that the weak factorization system $(\C_1\cap\W_1,\F_1)$ is small but we are not aware of a proof of this fact. Thus it is unlikely that the mixed model category $\top_m$ is combinatorial.

\bpf[Sketch of proof] It suffices to check that the factorization of a map by a strong cofibration which is a homotopy equivalence followed by a Hurewicz fibration is accessible. We can use the construction of \cite[Definition~3.2]{Barthel-Riel}. The middle space is given by an accessible functor as soon as the underlying category is locally presentable and cartesian closed. The result follows from Proposition~\ref{cole-model-structure}. 
\epf

We want to recall the following theorem: 

\bth (C. Rezk) \label{rezk-surprising-result} Let $(\C_1,\W_1,\F_1)$ and $(\C_2,\W_2,\F_2)$ be two model structures on the same underlying category with $\W_1=\W_2$. Then the model structure $(\C_1,\W_1,\F_1)$ is left proper (right proper resp.) if and only if the model structure $(\C_2,\W_2,\F_2)$ is left proper (right proper resp.).
\eth

\bpf This amazing result is a consequence of \cite[Proposition~2.5]{REZK200265} (it is \cite[Proposition~2.7]{REZK200265} in the prepublished version). 
\epf

It implies that $\top_m$ is proper. Indeed, it has the same class of weak equivalences as the model category $\top_Q$. And the latter is known to be proper by \cite[Theorem~13.1.11]{ref_model2}. Thus the former is proper as well by Theorem~\ref{rezk-surprising-result}. The mixed model structure $\top_m$ is also closed monoidal for the binary product by \cite[Proposition~6.6]{mixed-cole}.

\section{Enriched diagrams over a small enriched category}
\label{enricheddiag}

Let $\mathcal{P}$ be a nonempty enriched small category. Denote by $\mathcal{P}(\ell,\ell')$ the space of maps from $\ell$ to $\ell'$. The underlying category is denoted by $\mathcal{P}_0$ and we have \[\mathcal{P}_0(\ell,\ell')=\top(\{0\},\mathcal{P}(\ell,\ell'))\] for all objects $\ell$ and $\ell'$ of $\mathcal{P}$.

An enriched functor from $\mathcal{P}$ to $\top$ is a functor $F$ of $\dgr{\mathcal{P}}$ such that the set map  \[\mathcal{P}_0(\ell_1,\ell_2) \longrightarrow \top(F(\ell_1),F(\ell_2))\] induces a continuous map  \[F_{\ell_1,\ell_2}:\mathcal{P}(\ell_1,\ell_2) \longrightarrow \ttop(F(\ell_1),F(\ell_2)).\] An enriched natural transformation $\eta:F\Rightarrow G$ from an enriched functor $F$ to an enriched functor $G$ is, by definition \cite[Diagram~6.13]{Borceux2}, a family of continuous maps \[\eta_{\ell}:\{0\} \to \ttop(F(\ell),G(\ell))\] such that the following diagram of $\top$ commutes for all $\ell_1,\ell_2\in\Obj(\mathcal{P})$: 
	{\tiny\[
		\xymatrix
		{
			\mathcal{P}(\ell_1,\ell_2)\fr{(G_{\ell_1,\ell_2},\eta_{\ell_1})} \fd{(\eta_{\ell_2},F_{\ell_1,\ell_2})} & \ttop(G(\ell_1),G(\ell_2)) \p \ttop(F(\ell_1),G(\ell_1))\fd{}\\
			\ttop(F(\ell_2),G(\ell_2)) \p \ttop(F(\ell_1),F(\ell_2))\fr{} & \ttop(F(\ell_1),G(\ell_2))
		}
		\]}
Since $\top$ is cartesian closed, we have 
\[\top(F(\ell),G(\ell))) = \top(\{0\},\ttop(F(\ell),G(\ell))).\]
Therefore $\eta$ is just an ordinary natural transformation from $F$ to $G$ in $\dgr{\mathcal{P}}$. The underlying category $[\mathcal{P},\top]_0$ of the enriched category of enriched functors $[\mathcal{P},\top]$ can then be identified with a full subcategory of the category  $\dgr{\mathcal{P}}$ of functors $F:\mathcal{P}\to \top$ such that the set map  $\mathcal{P}_0(\ell_1,\ell_2) \longrightarrow \top(F(\ell_1),F(\ell_2))$ induces a continuous map  $\mathcal{P}(\ell_1,\ell_2) \longrightarrow \ttop(F(\ell_1),F(\ell_2))$ for all $\ell_1,\ell_2\in\Obj(\mathcal{P})$. 

It is well-known that the enriched category $[\mathcal{P},\top]$ is tensored and cotensored over $\top$ (e.g. see \cite[Lemma~5.2]{MoserLyne}). For an enriched diagram $F:\mathcal{P}\to \top$, and a topological space $U$, the enriched diagram $F\ot U:\mathcal{P}\to \top$ is defined by $F\ot U = F(-)\p U$ and $F^U:\mathcal{P}\to \top$ is defined by $F^U=\ttop(U,F(-))$.

\bp \label{enriched-loc} The category $[\mathcal{P},\top]_0$ is locally presentable. \ep

\bpf By Corollary~\ref{localbase}, the category $\top$ is a locally presentable base in the sense of \cite[Definition~1.1]{EnrichedSketch} for the closed monoidal structure given by the binary product. Using \cite[Example~6.2]{EnrichedSketch}, we deduce that the enriched category $[\mathcal{P},\top]$ is enriched locally presentable. The proof is complete using \cite[Proposition~6.6]{EnrichedSketch}.
\epf

\bp \label{car-enriched-functor} A functor $F$ of $\dgr{\mathcal{P}}$ belongs to the full subcategory $[\mathcal{P},\top]_0$ if and only if for all $\ell,\ell'\in\Obj(\mathcal{P})$, the set map $\mathcal{P}(\ell,\ell')\p F(\ell) \to F(\ell')$ defined by the mapping $(\phi,\gamma)\mapsto F(\phi)(\gamma)$ is continuous. \ep

\bpf This comes from the bijection of sets 
\[\top(\mathcal{P}(\ell,\ell'),\ttop(F(\ell),F(\ell')))\iso \top(\mathcal{P}(\ell,\ell') \p F(\ell),F(\ell')).\] \epf

Note that $\Delta_{\mathcal{P}}\varnothing$ belongs to $[\mathcal{P},\top]_0$ just because $\id_\varnothing$ is continuous.

\bp \label{colim-lim-enriched} The inclusion functor $[\mathcal{P},\top]_0 \subset \dgr{\mathcal{P}}$ is colimit-preserving and limit-preserving. \ep

\bpf
Since the category $[\mathcal{P},\top]_0$ is a full subcategory of $\dgr{\mathcal{P}}$, it suffices to prove that
$[\mathcal{P},\top]_0$  is closed under the colimits and the limits of $\dgr{\mathcal{P}}$. Let $(F_i)_{i\in I}$ be a small diagram of functors of $[\mathcal{P},\top]_0$. The case of colimits comes from the fact that the colimit of the maps 
\[\mathcal{P}(\ell,\ell')\p F_i(\ell) \to F_i(\ell')\]
in the category of diagrams $\Mor(\top)$ is 
\[\mathcal{P}(\ell,\ell')\p (\liminj F_i(\ell)) \to \liminj F_i(\ell')\]
because $\top$ is cartesian closed and because colimits in $\Mor(\top)$ are calculated pointwise. 
The case of limits comes from the fact that the limit of the maps 
\[\mathcal{P}(\ell,\ell')\p F_i(\ell) \to F_i(\ell')\]
in the category of diagrams $\Mor(\top)$ is 
\[\mathcal{P}(\ell,\ell')\p \limproj F_i(\ell) \to \limproj F_i(\ell')\]
because the functor $\limproj$ commutes with binary products as any right adjoint and because limits in $\Mor(\top)$ are calculated pointwise.
\epf

\begin{nota}
Let $\mathbb{F}^{\mathcal{P}}_{\ell}X=\mathcal{P}(\ell,-)\p X \in [\mathcal{P},\top]_0$ where $X$ is a topological space and where $\ell$ is an object of $\mathcal{P}$.
\end{nota}

\bp\label{ev-adj}
For every enriched functor $F:\mathcal{P}\to \top$, every $\ell\in\Obj(\mathcal{P})$ and every topological space $X$, we have 
the natural bijection of sets 
\[[\mathcal{P},\top]_0(\mathbb{F}^{\mathcal{P}}_{\ell}X,F) \iso \top(X,F(\ell)).\]
In particular, the functor $\mathbb{F}^{\mathcal{P}}_{\ell}:\top \to [\mathcal{P},\top]_0$ is colimit-preserving for all $\ell\in\Obj(\mathcal{P})$.
\ep

\bpf
We have the sequence of natural homeomorphisms
\begin{align}
[\mathcal{P},\top](\mathbb{F}^{\mathcal{P}}_{\ell}X,F) & \iso [\mathcal{P},\top](\mathcal{P}(\ell,-),\ttop(X,F(-)) \label{eq-1-1}\\&\iso \ttop(X,F(\ell))\label{eq-1-2},
\end{align}
$(\ref{eq-1-1})$ because the enriched category $[\mathcal{P},\top]$ is tensored and cotensored over $\top$, $(\ref{eq-1-2})$ by the enriched Yoneda lemma. By applying to the obtained isomorphism the functor $\top(\{0\},-)$, we obtain the desired bijection.
\epf

\begin{cor} \label{LLPdiag}
Let $f:X\to Y$ be a map of $\top$. The map of enriched diagrams $\mathbb{F}^{\mathcal{P}}_{\ell}X \to \mathbb{F}^{\mathcal{P}}_{\ell}Y$ induced by $f$ satisfies the LLP with respect to a map of diagrams $D\to E$ of $[\mathcal{P},\top]_0$ if and only if $f$ satisfies the LLP with respect to the continuous map $D_{\ell}\to E_{\ell}$. 
\end{cor}

\bth \label{cont-disc}
The inclusion functor $i^{\mathcal{P}}:[\mathcal{P},\top]_0 \subset \dgr{\mathcal{P}}$ has both a left adjoint and a right adjoint. In other terms, the category $[\mathcal{P},\top]_0$ is both a reflective and a coreflective subcategory of $\dgr{\mathcal{P}}$. If $i^{\mathcal{P}}_!:\dgr{\mathcal{P}}\to [\mathcal{P},\top]_0$ is the left adjoint, then for all $\ell\in\Obj(\mathcal{P})$ and all topological spaces $U$, $i^{\mathcal{P}}_!(\mathcal{P}_0(\ell,-)\p U) = \mathcal{P}(\ell,-)\p U$.
\eth

\bpf
Since the inclusion functor is colimit-preserving, it is in particular accessible and it is also a left adjoint because both the categories $[\mathcal{P},\top]_0$ and $\dgr{\mathcal{P}}$ are locally presentable. Since it is moreover limit-preserving, it is a right adjoint. We have the sequence of bijections 
\begin{align}
\label{cont-disc1}[\mathcal{P},\top]_0(i^{\mathcal{P}}_!(\mathcal{P}_0(\ell,-)\p U),Y)&\iso \dgr{\mathcal{P}}(\mathcal{P}_0(\ell,-)\p U,Y) \\
&\iso \int_{\ell'} \top(\mathcal{P}_0(\ell,\ell')\p U,Y(\ell'))\label{cont-disc2}\\
&\iso \int_{\ell'} \set(\mathcal{P}_0(\ell,\ell'),\top(U,Y(\ell')))\label{cont-disc3}\\
&\iso \set^{\mathcal{P}_0}(\mathcal{P}_0(\ell,-),\top(U,Y(-)))\label{cont-disc4}\\
&\iso \top(U,Y(\ell))\label{cont-disc5}\\
&\iso [\mathcal{P},\top]_0(\mathcal{P}(\ell,-)\p U,Y),\label{cont-disc6}
\end{align}
$(\ref{cont-disc1})$ because $[\mathcal{P},\top]_0$ is a full subcategory of $\dgr{\mathcal{P}}$ by Proposition~\ref{car-enriched-functor} and by adjunction, $(\ref{cont-disc2})$ by \cite[page 219 (2)]{MR1712872}, $(\ref{cont-disc3})$ because $\mathcal{P}_0(\ell,\ell')$ is a set, $(\ref{cont-disc4})$ by \cite[page 219 (2)]{MR1712872}, $(\ref{cont-disc5})$ by Yoneda, and finally $(\ref{cont-disc6})$ by Proposition~\ref{ev-adj}. The proof is complete thanks to the Yoneda lemma.
\epf

\bp \label{Fcolim}
Let $\ell\in\Obj(\mathcal{P})$. Let $U$ be a topological space. Then there is the natural homeomorphism  $\liminj_{\mathcal{P}} \mathbb{F}^{\mathcal{P}}_{\ell}U \iso U$.
\ep

\bpf There is the sequence of bijections ($V$ being any topological space): 
\begin{align}
\top(\liminj \mathbb{F}^{\mathcal{P}}_{\ell}U,&V) \iso \dgr{\mathcal{P}}(\mathbb{F}^{\mathcal{P}}_{\ell}U,\Delta_{\mathcal{P}}V)\label{l1}\\
&\iso [\mathcal{P},\top]_0(\mathbb{F}^{\mathcal{P}}_{\ell}U,\Delta_{\mathcal{P}}V)\label{l2}\\
&\iso \top\left(\{0\},[\mathcal{P},\top](\mathcal{P}(\ell,-)\p U,\Delta_{\mathcal{P}}V)\right)\label{l3}\\
&\iso \top\left(\{0\},[\mathcal{P},\top](\mathcal{P}(\ell,-),\ttop(U,\Delta_{\mathcal{P}}V(-)))\right)\label{l4}\\
&\iso \top\left(\{0\},\ttop(U,V)\right)\label{l5}\\
&\iso \top(U,V), \label{l6}
\end{align}
$(\ref{l1})$ by definition of the colimit, $(\ref{l2})$ because $[\mathcal{P},\top]_0$ is a full subcategory of $\dgr{\mathcal{P}}$ and because the constant diagram functor belongs to $[\mathcal{P},\top]_0$, $(\ref{l3})$ by definition of the enriched category $[\mathcal{P},\top]$, $(\ref{l4})$ since the enriched category $[\mathcal{P},\top]$ is tensored and cotensored over $\top$, $(\ref{l5})$ by the enriched Yoneda lemma, and finally $(\ref{l6})$ by definition of the enrichment of $\top$. The proof is complete thanks to the (ordinary) Yoneda lemma.
\epf

\section{The homotopy theory of enriched diagrams of topological spaces}
\label{homotopytheory}

\begin{nota}
Let $n\geq 1$. Denote by $\mathbf{D}^n = \{b\in \mathbb{R}^n, |b| \leq 1\}$ the $n$-dimensional disk, and by $\mathbf{S}^{n-1} = \{b\in \mathbb{R}^n, |b| = 1\}$ the $(n-1)$-dimensional sphere. By convention, let $\mathbf{D}^{0}=\{0\}$ and $\mathbf{S}^{-1}=\varnothing$. 
\end{nota}

\bth \label{proj-model-structure-enriched-diagram} The category $[\mathcal{P},\top]_0$ can be endowed with a structure of a combinatorial model category as follows: 
\begin{itemize}
\item The set of generating cofibrations is the set of maps \[\{\mathbb{F}^{\mathcal{P}}_{\ell}\mathbf{S}^{n-1}\to \mathbb{F}^{\mathcal{P}}_{\ell}\mathbf{D}^n\mid n\geq 0,\ell\in\Obj(\mathcal{P})\}\]
induced by the inclusions $\mathbf{S}^{n-1}\subset \mathbf{D}^n$.
\item The set of generating trivial cofibrations is the set of maps \[\{\mathbb{F}^{\mathcal{P}}_{\ell}\mathbf{D}^{n}\to \mathbb{F}^{\mathcal{P}}_{\ell}\mathbf{D}^{n+1}\mid n\geq 0,\ell\in\Obj(\mathcal{P})\}\]
where the maps $\mathbf{D}^{n}\subset \mathbf{D}^{n+1}$ are induced by the mappings \[(x_1,\dots,x_n) \mapsto (x_1,\dots,x_n,0). \]
\item A map $F\to G$ is a weak equivalence if and only if for all $\ell\in\Obj(\mathcal{P})$, the continuous map $F(\ell)\to G(\ell)$ is a weak equivalence of $\top_Q$, i.e. the weak equivalences are the pointwise weak homotopy equivalences
\item A map $F\to G$ is a fibration if and only if for all $\ell\in\Obj(\mathcal{P})$, the continuous map $F(\ell)\to G(\ell)$ is a fibration of $\top_Q$, i.e. the fibrations are the pointwise Serre fibrations.
\end{itemize}
This model structure, denoted by $[\mathcal{P},\top_Q]^{proj}_0$, is called the {\rm projective mo\-del structure}. The cofibrations are called the projective cofibrations.
\eth

\bpf The existence of an accessible model structure is a consequence of Proposition~\ref{enriched-loc}, \cite[Theorem~6.5(ii)]{MoserLyne} and of the fact that all objects of $\top_Q$ are fibrant (it suffices to use the adjunction $\top^{\Obj(\mathcal{P})}\leftrightarrows [\mathcal{P},\top]_0$ and the Quillen Path Object Argument, which implies the acyclicity condition). It is cofibrantly generated by Corollary~\ref{LLPdiag} and because the set of inclusions $\{\mathbf{S}^{n-1}\subset \mathbf{D}^n\mid n\geq 0\}$ ($\{\mathbf{D}^{n}\subset \mathbf{D}^{n+1}\mid n\geq 0\}$ resp.) is a set of generating (trivial resp.) cofibrations of $\top_Q$. 
\epf

In \cite{PiaDgr}, Piacenza proves a similar result by working in the category of Hausdorff $k$-spaces in the sense of \cite{MR45:9323}. We do not know if Piacenza's proof can be adapted to $\Delta$-generated spaces, especially because Piacenza works with Hausdorff spaces and we do not assume any separation hypothesis. It is also known that Hausdorff $k$-spaces do not behave very well for algebraic topology problems and that weak Hausdorff $k$-spaces are much better (see the end of the introduction of \cite{mdtop} for some bibliographical research about this problem).

\cite[Theorem 24.4]{shuhomotopycolim} and \cite[Theorem~4.4]{MoserLyne} (the latter is a generalization of the former to the framework of accessible model categories) give sufficient conditions for the projective model structure to exist in an enriched situation. They could be applied to prove Theorem~\ref{proj-model-structure-enriched-diagram} if e.g. we assumed that all topological spaces $\mathcal{P}(\ell,\ell')$ were cofibrant in $\top_Q$. 

\begin{cor} \label{uneadj} The adjunction $i^{\mathcal{P}}_!\dashv i^{\mathcal{P}}$ of Theorem~\ref{cont-disc} is a Quillen adjunction between the projective model structures of $\dgr{\mathcal{P}}$ and $[\mathcal{P},\top]_0$.
\end{cor}

\bth Suppose that all topological spaces $\mathcal{P}(\ell,\ell')$ are homotopy equivalent to a cofibrant space. The category $[\mathcal{P},\top_m]_0$ can be endowed with a structure of accessible model category characterized as follows: 
\begin{itemize}
\item A map $F\to G$ is a cofibration if and only if for all $\ell\in\Obj(\mathcal{P})$, the continuous map $F(\ell)\to G(\ell)$ is a cofibration of $\top_m$, i.e. the cofibrations are the pointwise mixed cofibrations.
\item A map $F\to G$ is a weak equivalence if and only if for all $\ell\in\Obj(\mathcal{P})$, the continuous map $F(\ell)\to G(\ell)$ is a weak equivalence of $\top_m$, i.e. the weak equivalences are the pointwise weak homotopy equivalences
\end{itemize}
This model structure, denoted by $[\mathcal{P},\top_m]^{inj}_0$, is called the {\rm injective mixed model structure}. The fibrations are called the injective mixed fibrations. 
\eth

\bpf
By \cite[Proposition~6.4]{mixed-cole}, the cartesian closed category $\top$ equi\-pped with the mixed model structure is a monoidal model category. The proof is complete by Proposition~\ref{enriched-loc}, Corollary~\ref{mixed-accessible} and \cite[Theorem~6.5(i)]{MoserLyne}
\epf

\begin{cor} \label{leftproper-first} Suppose that all topological spaces $\mathcal{P}(\ell,\ell')$ are homotopy equivalent to a cofibrant space. Then the projective model structure $[\mathcal{P},\top_Q]^{proj}_0$ is proper. \end{cor}

\bpf By hypothesis, all topological spaces $\mathcal{P}(\ell,\ell')$ are mixed cofibrant. We deduce that the functors $-\p \mathcal{P}(\ell,\ell'):\top_m\to\top_m$ preserve mixed cofibrations for all $\ell,\ell'\in  \Obj(\mathcal{P})$ because $(\top_m,\p)$ is a closed monoidal model category. Using \cite[Proposition~8.1(i)]{MoserLyne}, we obtain that the injective model structure $[\mathcal{P},\top_m]^{inj}_0$ is left proper. By Theorem~\ref{rezk-surprising-result}, we deduce that $[\mathcal{P},\top_Q]^{proj}_0$ is left proper\footnote{We cannot apply \cite[Proposition~8.1(i)]{MoserLyne} directly to $[\mathcal{P},\top_Q]^{proj}_0$ because the spaces of maps of $\mathcal{P}$ are not necessarily cofibrant in the Quillen model structure of $\top$.}. The model category $[\mathcal{P},\top_Q]^{proj}_0$ is right proper because the fibrations are the pointwise fibrations and because all topological spaces are fibrant.
\epf

We can actually remove the hypothesis of Corollary~\ref{leftproper-first} but the proof is a little bit more involved. It makes use of \cite[Proposition~A.2]{mdtop} and \cite[Proposition~A.6]{mdtop}; the proof of \cite[Proposition~A.6]{mdtop} is based on the notion of relative $T_1$-inclusion of \cite{hocolimfacile}. We include it here for completness.

\bth \label{lefproper-ok} The projective model structure $[\mathcal{P},\top_Q]^{proj}_0$ is proper. \eth

\bpf It suffices to prove that it is left proper because all objects are fibrant. In a model category, weak equivalences are closed under retract. Therefore it suffices to prove that the pushout of a weak equivalence along a transfinite composition of pushouts of maps of the form $\mathbb{F}^{\mathcal{P}}_\ell \mathbf{S}^{n-1}\to \mathbb{F}^{\mathcal{P}}_\ell \mathbf{D}^{n}$ is still a weak equivalence. Consider first the following situation: 
\[
\xymatrix
{
\mathbb{F}^{\mathcal{P}}_\ell \mathbf{S}^{n-1} \fr{} \fd{} & F \fr{f}\fd{} & G \fd{} \\
\mathbb{F}^{\mathcal{P}}_\ell \mathbf{D}^{n} \fr{} & \cocartesien H \fr{\widetilde{f}}& \cocartesien K.
}
\]
For all objects $\ell'$ of $\mathcal{P}$, we obtain the diagram of $\top$
\[
\xymatrix
{
\mathcal{P}(\ell,\ell')\p \mathbf{S}^{n-1} \fr{} \fd{} & F(\ell') \fr{f_{\ell'}}\fd{} & G(\ell') \fd{} \\
\mathcal{P}(\ell,\ell')\p \mathbf{D}^{n} \fr{} & \cocartesien H(\ell') \fr{\widetilde{f}_{\ell'}}& \cocartesien K(\ell').
}
\]
If $f_{\ell'}$ is a weak homotopy equivalence, then $\widetilde{f}_{\ell'}$ is a weak homotopy equivalence by \cite[Proposition~A.2]{mdtop}. Thus if $f$ is a pointwise weak equivalence, then $\widetilde{f}$ is a pointwise weak equivalence. By \cite[Proposition~A.6]{mdtop}, this process can be iterated transfinitely since colimits in $[\mathcal{P},\top]_0$ are calculated pointwise by Proposition~\ref{colim-lim-enriched}. 
\epf

\bp \label{example-proj-cof}
For all (trivial resp.) cofibrations $f:U\to V$ of $\top_Q$ and all $\ell\in\Obj(\mathcal{P})$, the map of diagrams $\mathbb{F}^{\mathcal{P}}_{\ell}U \to \mathbb{F}^{\mathcal{P}}_{\ell}V$ is a (trivial resp.) projective cofibration. 
\ep

\bpf 
The model structure on $[\mathcal{P},\top]_0$ is obtained by right-inducing it using the adjunction $\top^{\Obj(\mathcal{P})}\leftrightarrows [\mathcal{P},\top]_0$ (see the comment in the proof of Theorem~\ref{proj-model-structure-enriched-diagram}) where the right adjoint is given by the evaluation maps. The proof is complete thanks to Proposition~\ref{ev-adj}. 
\epf

\begin{cor} \label{tractable}
The combinatorial model category $[\mathcal{P},\top_Q]^{proj}_0$ is tractable (i.e. the generating cofibrations and generating trivial cofibrations have cofibrant domains).
\end{cor}

\bpf It is a consequence of the fact that the maps $\varnothing\subset \mathbf{S}^{n-1}$ and $\varnothing\subset \mathbf{D}^{n}$ are cofibrations for all $n\geq 0$ and that  $\mathbb{F}^{\mathcal{P}}_{\ell}\varnothing = \Delta_{\mathcal{P}}\varnothing$ is the initial object of $[\mathcal{P},\top]_0$ for all $\ell\in\Obj(\mathcal{P})$. The proof is complete thanks to Proposition~\ref{example-proj-cof}. \epf

\section{The main theorem}
\label{mainthm}

\bp \label{proj-cof-pointwise-mixed} Suppose that all topological spaces $\mathcal{P}(\ell,\ell')$ are homotopy equivalent to a cofibrant space. Then any cofibration of $[\mathcal{P},\top_Q]^{proj}_0$ is a cofibration of $[\mathcal{P},\top_m]^{inj}_0$. In other terms, the identity functor is a left Quillen adjoint \[\id:[\mathcal{P},\top_Q]^{proj}_0 \to [\mathcal{P},\top_m]^{inj}_0.\]  \ep

Another way to formulate this proposition is that any projective cofibration of the model category $[\mathcal{P},\top_Q]^{proj}_0$ is a pointwise mixed cofibration of $\top_m^{\mathcal{P}_0}$. On the contrary, a projective cofibration of $[\mathcal{P},\top_Q]^{proj}_0$ is not necessarily a pointwise cofibration of $\dgrQ{\mathcal{P}}$. For example, the diagram $\mathbb{F}^{\mathcal{P}}_{\ell}U$ is projective cofibrant for any cofibrant space $U$ by Proposition~\ref{example-proj-cof}. But the vertices of this diagram are only homotopy equivalent to a cofibrant space. 

\bpf Every cofibration of $[\mathcal{P},\top_Q]^{proj}_0$ is a retract of a transfinite composition of pushouts of maps of \[\{\mathbb{F}^{\mathcal{P}}_{\ell}\mathbf{S}^{n-1}\subset \mathbb{F}^{\mathcal{P}}_{\ell}\mathbf{D}^{n}\mid n\geq 0,\ell\in\Obj(\mathcal{P})\}.\] Therefore it suffices to prove that the maps $\mathbb{F}^{\mathcal{P}}_{\ell}\mathbf{S}^{n-1}\subset \mathbb{F}^{\mathcal{P}}_{\ell}\mathbf{D}^{n}$ are pointwise mixed cofibrations of $\top^{\mathcal{P}_0}$ for all $n\geq0,\ell\in\Obj(\mathcal{P})$. It suffices to prove that for all $\ell,\ell'\in\Obj(\mathcal{P})$ and all $n\geq 0$, the map \[\mathcal{P}(\ell,\ell') \p \mathbf{S}^{n-1} \subset \mathcal{P}(\ell,\ell')\p \mathbf{D}^{n}\] is a mixed cofibration. The latter fact comes from the facts that $\mathcal{P}(\ell,\ell')$ is cofibrant in $\top_m$, that any cofibration is a mixed cofibration and that $(\top_m,\p)$ is a monoidal model structure.
\epf

\begin{cor}
Suppose that all topological spaces $\mathcal{P}(\ell,\ell')$ are homotopy equivalent to a cofibrant space. Then the identity functor induces a left Quillen equivalence \[\id:[\mathcal{P},\top_Q]^{proj}_0 \to [\mathcal{P},\top_m]^{inj}_0.\]
\end{cor}

\bp \label{inj} There is a Quillen adjunction $\liminj\dashv \Delta_{\mathcal{P}}$ between the model categories $[\mathcal{P},\top_Q]^{proj}_0$ and $\top_Q$. 
\ep

\bpf
There is the sequence of bijections ($X$ being an object of $[\mathcal{P},\top]_0$ and $U$ being a topological space)
\[\top(\liminj X,U) \iso \dgr{\mathcal{P}}(X,\Delta_{\mathcal{P}}U) 
\iso [\mathcal{P},\top]_0(X,\Delta_{\mathcal{P}}U),\]
the left-hand bijection by definition of the colimit and by Proposition~\ref{colim-lim-enriched}, the right-hand bijection because the category $[\mathcal{P},\top]_0$ is a full subcategory of $\dgr{\mathcal{P}}$ and because the constant diagram functor belongs to $[\mathcal{P},\top]_0$. The right adjoint $\Delta_{\mathcal{P}}:\top_Q\to [\mathcal{P},\top_Q]^{proj}_0$ takes (trivial resp.) fibrations to pointwise (trivial resp.) fibrations, therefore it is a right Quillen adjoint.
\epf

\bp \label{inj2} There is a categorical adjunction $\Delta_{\mathcal{P}} \dashv \limproj$ between the categories $[\mathcal{P},\top]_0$ and $\top$. 
\ep

\bpf
The proof is analogous to the first part of the proof of Proposition~\ref{inj}.
\epf

\begin{cor}
Suppose that all topological spaces $\mathcal{P}(\ell,\ell')$ are homotopy equivalent to a cofibrant space. There is a Quillen adjunction $\Delta_{\mathcal{P}} \dashv \limproj$ between the model categories $\top_m$ and $[\mathcal{P},\top_m]^{inj}_0$.
\end{cor}

\bpf
The left adjoint $\Delta_{\mathcal{P}}:\top_m\to [\mathcal{P},\top_m]^{inj}_0$ takes (trivial resp.) mixed cofibrations to pointwise (trivial resp.) mixed cofibrations. We deduce that it is a left Quillen adjoint.
\epf

We can now prove the main theorem of the paper.

\bth \label{eq-topdgr-top} Suppose that all spaces $\mathcal{P}(\ell,\ell')$ are contractible. Then the Quillen adjunction 
\[\adj{[\mathcal{P},\top_Q]^{proj}_0}{\liminj}{\Delta_{\mathcal{P}}}{\top_Q}{3}{20}\]
is a Quillen equivalence.
\eth



\bpf Since all spaces $\mathcal{P}(\ell,\ell')$ are contractible by hypothesis, i.e. homotopy equivalent to a point, they are cofibrant for the mixed model structure $\top_m$. Let $U$ be a topological space. Let $U^{cof}\to U$ be a cofibrant replacement of $U$ in $\top_Q$. We obtain a map \[U^{cof} \longrightarrow (\Delta_{\mathcal{P}}U)(\ell)\] for some $\ell\in\Obj(\mathcal{P})$. By Proposition~\ref{ev-adj}, we obtain a map \[\mathbb{F}^{\mathcal{P}}_{\ell}U^{cof}\longrightarrow \Delta_{\mathcal{P}}U\] of $[\mathcal{P},\top]_0$. Since all topological spaces $\mathcal{P}(\ell,\ell')$ are contractible by hypothesis, the map \[\mathbb{F}^{\mathcal{P}}_{\ell}U^{cof}\longrightarrow \Delta_{\mathcal{P}}U\] is a weak equivalence of $[\mathcal{P},\top_Q]^{proj}_0$. By Proposition~\ref{example-proj-cof}, the map \[\mathbb{F}^{\mathcal{P}}_{\ell}\varnothing\to\mathbb{F}^{\mathcal{P}}_{\ell}U^{cof}\] is a cofibration of $[\mathcal{P},\top_Q]^{proj}_0$. Thus $\mathbb{F}^{\mathcal{P}}_{\ell}U^{cof}$ is a cofibrant replacement of $\Delta_{\mathcal{P}}U$ in $[\mathcal{P},\top_Q]^{proj}_0$. Using Proposition~\ref{Fcolim}, we obtain that the canonical map \[\liminj \mathbb{F}^{\mathcal{P}}_{\ell}U^{cof} \iso U^{cof} \longrightarrow U\] is a weak equivalence of $\top_Q$. We deduce that the functor \[\liminj:[\mathcal{P},\top_Q]^{proj}_0\longrightarrow\top_Q\] is homotopically surjective. 

Let $Y$ be a cofibrant diagram of $[\mathcal{P},\top_Q]^{proj}_0$. Since all objects of $[\mathcal{P},\top_Q]^{proj}_0$ are fibrant, we need to prove that the unit of the adjunction \[Y \longrightarrow \Delta_{\mathcal{P}}(\liminj Y)\] is a weak equivalence to complete the proof. Every cofibrant diagram of $[\mathcal{P},\top_Q]^{proj}_0$ is a retract of a cellular object of $[\mathcal{P},\top_Q]^{proj}_0$, i.e. of a transfinite composition of pushouts of generating cofibrations. And in a model category, the retract of a weak equivalence is a weak equivalence. We can therefore assume without lack of generality that $Y$ is cellular. As a first step, consider the commutative diagram in $[\mathcal{P},\top]_0$
\[
\xymatrix
{
	\mathbb{F}^{\mathcal{P}}_{\ell}\mathbf{S}^{n-1} \ar@{->}[d]\ar@{->}[ddr]\ar@{->}[r] & X \ar@{->}[ddr]  \\
	\mathbb{F}^{\mathcal{P}}_{\ell}\mathbf{D}^{n}\ar@{->}[ddr]  & \\
	& \Delta_{\mathcal{P}}\liminj \mathbb{F}^{\mathcal{P}}_{\ell}\mathbf{S}^{n-1}\ar@{->}[d] \ar@{->}[r] & \Delta_{\mathcal{P}}\liminj X \\
	& \Delta_{\mathcal{P}}\liminj \mathbb{F}^{\mathcal{P}}_{\ell}\mathbf{D}^{n}  & 
}
\]
obtained using the unit of the adjunction $\id \Rightarrow \Delta_{\mathcal{P}}\liminj$. Suppose that the map \[X\longrightarrow \Delta_{\mathcal{P}}\liminj X\] is a pointwise weak equivalence (of the model category $[\mathcal{P},\top_Q]^{proj}_0$ or equivalently of the model category $[\mathcal{P},\top_m]^{inj}_0$), that $X$ is projective cofibrant and that $\Delta_{\mathcal{P}}\liminj X$ is pointwise mixed cofibrant. The map \[\mathbb{F}^{\mathcal{P}}_{\ell}\mathbf{S}^{n-1}\to \mathbb{F}^{\mathcal{P}}_{\ell}\mathbf{D}^{n}\] is a pointwise mixed cofibration between pointwise mixed cofibrant diagrams by Proposition~\ref{proj-cof-pointwise-mixed}. By Proposition~\ref{Fcolim}, the map \[\Delta_{\mathcal{P}}\liminj \mathbb{F}^{\mathcal{P}}_{\ell}\mathbf{S}^{n-1}\iso \Delta_{\mathcal{P}}\mathbf{S}^{n-1}\to \Delta_{\mathcal{P}}\mathbf{D}^{n}\iso \Delta_{\mathcal{P}}\liminj \mathbb{F}^{\mathcal{P}}_{\ell}\mathbf{D}^{n}\] is a pointwise mixed cofibration between pointwise mixed cofibrant diagrams as well because all Quillen cofibrations are mixed cofibrations. Since $X$ is projective cofibrant by hypothesis, it is also pointwise mixed cofibrant by Proposition~\ref{proj-cof-pointwise-mixed}.  For all topological spaces $U$, the maps \[\mathbb{F}^{\mathcal{P}}_{\ell}U\to \Delta_{\mathcal{P}}U\] are pointwise weak equivalences since all spaces $\mathcal{P}(\ell,\ell')$ are contractible. We are ready to apply the cube lemma \cite[Lemma~5.2.6]{MR99h:55031} in $[\mathcal{P},\top_m]^{inj}_0$ by passing to the colimit. Since the functor \[\Delta_{\mathcal{P}}.\liminj:[\mathcal{P},\top]_0\to[\mathcal{P},\top]_0\] is colimit-preserving by Proposition~\ref{inj} and Proposition~\ref{inj2} as a composite of two colimit-preserving functors, we obtain the commutative diagram
\[
\xymatrix
{
	\mathbb{F}^{\mathcal{P}}_{\ell}\mathbf{S}^{n-1} \ar@{->}[d]\ar@{->}[ddr]\ar@{->}[r] & X \ar@{->}[ddr] \ar@{->}[d] \\
	\mathbb{F}^{\mathcal{P}}_{\ell}\mathbf{D}^{n}\ar@{->}[ddr] \ar@{->}[r] & Y \cocartesien \ar@{-->}[ddr]\\
	& \Delta_{\mathcal{P}}\mathbf{S}^{n-1}\ar@{->}[d] \ar@{->}[r] & \Delta_{\mathcal{P}}\liminj X \ar@{->}[d]\\
	& \Delta_{\mathcal{P}}\mathbf{D}^{n} \ar@{->}[r] & \cocartesien \Delta_{\mathcal{P}}\liminj Y.
}
\]
Using the cube lemma in $[\mathcal{P},\top_m]^{inj}_0$, we deduce that the map \[Y\longrightarrow \Delta_{\mathcal{P}}\liminj Y\] is a pointwise weak equivalence from a projective cofibrant object of $[\mathcal{P},\top_Q]^{proj}_0$ to a pointwise mixed cofibrant object of $[\mathcal{P},\top_m]^{inj}_0$. Moreover, the map $X\to Y$ is a cofibration both of $[\mathcal{P},\top_Q]^{proj}_0$ and of $[\mathcal{P},\top_m]^{inj}_0$, and the map \[\Delta_{\mathcal{P}}\liminj X\to \Delta_{\mathcal{P}}\liminj Y\] is a pointwise mixed cofibration, i.e. a cofibration of $[\mathcal{P},\top_m]^{inj}_0$ as well. We start from $X=\varnothing$: the only possibility is then $n=0$, $\mathbb{F}^{\mathcal{P}}_{\ell}\mathbf{S}^{-1}=\varnothing$ is the initial enriched diagram and $Y=\mathbb{F}^{\mathcal{P}}_{\ell}\mathbf{D}^{0}$ is a free enriched diagram generated by a point (it depends on $\ell$). Then by iterating the process transfinitely, we obtain two transfinite towers of cofibrations of $[\mathcal{P},\top_m]^{inj}_0$ between pointwise mixed cofibrant diagrams. In this case, the colimit is a homotopy colimit. Thus, for all cellular objects $Y$ of $[\mathcal{P},\top_Q]^{proj}_0$, the unit of the adjunction \[Y\longrightarrow \Delta_{\mathcal{P}}\liminj Y\] is a pointwise weak homotopy equivalence. 
\epf

\section{Concluding remarks}
\label{conclusion}

We conclude this paper with two remarks.

\subsection*{About the monoid of nondecreasing continuous maps from $[0,1]$ to itself}
The behavior of the homotopy colimit $\holiminj \D^{\mathcal{M}}(X,dX)$ in the enriched case remains the same because Theorem~\ref{eq-topdgr-top} can still be applied. However, the homotopy colimit \textit{in the non-enriched setting} behaves well in this case. To see that, we start from the surprising observation due to Tyler Lawson (let us recall that unlike groups, every Quillen cofibrant path-connected space has the same homotopy type as the classifying space of some monoid \cite[Theorem~1]{McDuff}): 

%

\bth \cite{311786} (T. Lawson) \label{bm-contractible}
The classifying space $B\mathcal{M}$ is contractible.
\eth

Note that in the proof below, the two maps of monoids $U,V:\mathcal{M}\to \mathcal{M}$ induce two endomorphisms of the group $\mathcal{G}$. However, $\phi$ and $\psi$, which provide the homotopies, do not belong to $\mathcal{G}$. 

\bpf
Consider $\phi,\psi\in \mathcal{M}$ defined by 
\begin{align*}
	\phi(x) =  &\begin{cases}2x &\hbox{if }x\leq 1/2 \\
	1 &\hbox{if }x\geq 1/2
	\end{cases}\\
	\psi(x) =  &\begin{cases}0 &\hbox{if }x\leq 1/2 \\
	2x-1 &\hbox{if }x\geq 1/2
	\end{cases}
\end{align*}
Define two maps of monoids $U,V:\mathcal{M}\to \mathcal{M}$ by 
\begin{align*}
(Uf)(x) &=
\begin{cases} \tfrac{1}{2}f(2x) &\text{if }x \leq 1/2\\x &\text{if }x \geq 1/2\end{cases}\\
(Vf) &= \id_{[0,1]}
\end{align*}
For any $f\in \mathcal{M}$, we have the following identities: 
\begin{align*}
\phi . (Uf) &= \id_{\mathcal{M}}(f) . \phi\\
\psi . (Uf) &= (Vf) . \psi
\end{align*}
As a result, we can reinterpret this in terms of the one-object category $\mathcal{M}$: we get three functors $\id_{\mathcal{M}},U,V: \mathcal{M} \to \mathcal{M}$ and natural transformations $\phi: U \Rightarrow \id_{\mathcal{M}}$ and $\psi: U \Rightarrow V$. Upon taking geometric realization, we get a space $B\mathcal{M}$, these functors turn into continuous maps $\id_{B\mathcal{M}}, BU, BV:B\mathcal{M} \to B\mathcal{M}$ and homotopies $B\phi$ from $BU$ to $\id_{B\mathcal{M}}$ and $B\psi$ from $BU$ to $BV$. However, $BV$ is a constant map, and so this says that the homotopy type of $B\mathcal{M}$ is contractible.
\epf

We obtain as an immediate consequence: 

\bth \label{good-behavior} We have a weak homotopy equivalence \[\holiminj \D^{\mathcal{M}}(X,dX)\simeq dX\] where the homotopy colimit is calculated in the ordinary (i.e. non-enriched) projective model structure. \eth

\bpf
Every map $\phi^*:dX\to dX$ is homotopic to the identity by the homotopy $H:[0,1]\p dX\to dX$ taking $(t,\gamma)$ to $\gamma.(t.\phi+(1-t).\id_{[0,1]})$. We deduce that every map of $\D^{\mathcal{M}}(X,dX)$ is a weak homotopy equivalence. Since $B\mathcal{M}$ is contractible by Theorem~\ref{bm-contractible}, the ordinary (i.e. non-enriched) homotopy colimit $\holiminj \D^{\mathcal{M}}(X,dX)$ is in this case weakly homotopy equivalent to $dX$ by \cite[Corollary~29.2]{monographie_hocolim}.
\epf

\subsection*{Toward a generalization of the main theorem}
Every enriched functor $f:\mathcal{P}_1\to \mathcal{P}_2$ gives rise to a Quillen adjunction \[f_*\dashv f^*:[\mathcal{P}_1,\top_Q]^{proj}_0 \leftrightarrows [\mathcal{P}_2,\top_Q]^{proj}_0\] where $f^*(G)=G.f$ is the precomposition by $f$ and where \[f_*F(-)=\int^{\ell'} \mathcal{P}_2(f(\ell'),-)\p F(\ell')\] is the enriched left Kan extension. The main theorem can be reformulated as follows: 

\bth \label{gen}
Suppose that $\mathcal{P}$ is locally contractible. Let $f:\mathcal{P}\to \thin{\mathcal{P}}$ be the unique functor where $\thin{\mathcal{P}}$ is the enriched small category which has the same objects as $\mathcal{P}$ and exactly one map between each object. Then the functor \[f^*:[\thin{\mathcal{P}},\top_Q]^{proj}_0\to [\mathcal{P},\top_Q]^{proj}_0\] is a right Quillen equivalence.
\eth

\bpf
By Theorem~\ref{eq-topdgr-top}, there is the right Quillen equivalence \[\Delta_{\thin{\mathcal{P}}}:\top_Q\longrightarrow [\thin{\mathcal{P}},\top_Q]^{proj}_0.\] The composite functor \[\top_Q \to [\thin{\mathcal{P}},\top_Q]^{proj}_0 \stackrel{f^*}\to [\mathcal{P},\top_Q]_0\] is the constant diagram functor $\Delta_{\mathcal{P}}$ which is a right Quillen equivalence by Theorem~\ref{eq-topdgr-top}. Since $f^*$ preserves pointwise fibrations and pointwise weak equivalences, it is a right Quillen functor. By the $2$-out-of-$3$ property, we deduce that it is a right Quillen equivalence. 
\epf

This raises the question of generalizing Theorem~\ref{gen} by considering an enriched functor $f:\mathcal{P}_1\to \mathcal{P}_2$ between enriched small categories which is essentially surjective and locally a weak homotopy equivalence. Similar questions are studied in \cite[Proposition~22.5 and Proposition~22.9]{shuhomotopycolim} in the context of enriched homotopical categories. Unfortunately, Shulman's work cannot be used here, at least without adaptation, because the goodness condition is not satisfied. The main obstacle is that the injection of a point need not be a cofibration (cf. \cite[Definition~23.11]{shuhomotopycolim}).

\appendix

\section{Variant for the particular cases of \texorpdfstring{$\mathcal{G}$ and $\mathcal{M}$}{Lg}}

The obstacle mentioned in Section~\ref{conclusion} to apply \cite{shuhomotopycolim} does not exist in the particular case of $\mathcal{G}$ and $\mathcal{M}$. Indeed we have the following proposition: 

\bp \label{well-pointed}
For any map $\phi$ of  $\mathcal{G}$, the injection $\{\phi\}\subset \mathcal{G}$ is a mixed cofibration. For any map $\phi$ of  $\mathcal{M}$, the injection $\{\phi\}\subset \mathcal{M}$ is a mixed cofibration. 
\ep

\bpf We write the proof for $\mathcal{G}$. It is similar for $\mathcal{M}$. The homotopy $H:\mathcal{G}\p [0,1]\to \mathcal{G}$ defined by $(\psi,t)\mapsto t\psi +(1-t)\phi$ between $\phi$ and $\id_{\mathcal{G}}$ satisfies $H(\phi,t)=\phi$ for all $t\in[0,1]$. Thus the inclusion $\{\phi\}\subset \mathcal{G}$ is an inclusion of a strong deformation retract. Let $k_\Delta:\mathcal{TOP}\to \top$ be the right adjoint of the inclusion functor $\top\subset \mathcal{TOP}$ (the $\Delta$-kelleyfication functor). The space $\mathcal{G}$ is equipped with the $\Delta$-kelleyfication of the initial topology making the inclusion $\mathcal{G}\subset \ttop([0,1],[0,1])$ continuous. The space $\ttop([0,1],[0,1])$ is equal to $k_\Delta(\mathcal{TOP}([0,1],[0,1]))$ where $\mathcal{TOP}([0,1],[0,1])$ is the set of continuous maps from $[0,1]$ to $[0,1]$ equipped with the compact-open topology. We have a composite of continuous maps in $\mathcal{TOP}$ 
\[\mathcal{G} \longrightarrow \ttop([0,1],[0,1]) \longrightarrow \mathcal{TOP}([0,1],[0,1])\]
since the underlying set of $\ttop([0,1],[0,1])$ and $\mathcal{TOP}([0,1],[0,1])$ are the same and $\ttop([0,1],[0,1])$ have more open sets than the compact-open topology. Since $[0,1]$ is compact metrisable with the Euclidian metric defined by $d_2(x,y)=|x-y|$, the topological space  $\mathcal{TOP}([0,1],[0,1])$ is metrisable, and the compact-open topology is induced by the metric \[d(f,g)=\max\limits_{x\in [0,1]} d_2(f(x),g(x))\] by \cite[Proposition~A.13]{MR1867354}. We obtain a composite of continuous maps 
\[\xymatrix@C=2em{q:\mathcal{G} \fr{}& \ttop([0,1],[0,1]) \fr{}& \mathcal{TOP}([0,1],[0,1]) \fR{f\mapsto d(f,\phi)}&& [0,1]}\]
such that $q^{-1}(0)=\{\phi\}$. By \cite[Theorem~3]{MR0211403} (see also \cite[Proposition~2.3]{Barthel-Riel}), we deduce that the inclusion $\{\phi\}\subset \mathcal{G}$ is a cofibration of the Cole-Str{\o}m model structure. Since this map is homotopic to $\id_{\{\phi\}}$ which is a cofibration of $\top_Q$, we deduce by \cite[Proposition~3.6]{mixed-cole} that the map $\{\phi\}\subset \mathcal{G}$ is a mixed cofibration. 
\epf

\bth \label{proj-mixed-model-structure-enriched-diagram} Let $\mathcal{P}$ be an enriched small category. The category $[\mathcal{P},\top]_0$ can be endowed with a structure of accessible model category characterized as follows: 
\begin{itemize}
\item A map $F\to G$ is a weak equivalence if and only if for all $\ell\in\Obj(\mathcal{P})$, the continuous map $F(\ell)\to G(\ell)$ is a weak equivalence of $\top_m$, i.e. the weak equivalences are the pointwise weak homotopy equivalences
\item A map $F\to G$ is a fibration if and only if for all $\ell\in\Obj(\mathcal{P})$, the continuous map $F(\ell)\to G(\ell)$ is a fibration of $\top_m$, i.e. the fibrations are the pointwise Hurewicz fibrations.
\end{itemize}
This model structure, denoted by $[\mathcal{P},\top_m]^{proj}_0$, is called the {\rm projective mixed model structure}. The cofibrations are called the projective mixed cofibrations. This model structure is proper. It is Quillen equivalent to the projective model structure $[\mathcal{P},\top_Q]^{proj}_0$
\eth

\bpf The existence of an accessible model structure is a consequence of Proposition~\ref{enriched-loc}, \cite[Theorem~6.5(ii)]{MoserLyne} and of the fact that all objects of $\top_m$ are fibrant. It is proper by Theorem~\ref{rezk-surprising-result} and Theorem~\ref{lefproper-ok}. The identity functor induces a left Quillen adjoint \[[\mathcal{P},\top_Q]^{proj}_0 \to [\mathcal{P},\top_m]^{proj}_0\] since every projective (trivial) cofibration is a projective mixed (trivial) cofibration. All objects are fibrant and the two model categories have the same class of weak equivalences. Therefore, it is a Quillen equivalence by \cite[Definition~8.5.20]{ref_model2}. 
\epf

\bth \label{gen2}
Suppose that $\mathcal{P}$ is either $\mathcal{G}$ or $\mathcal{M}$ viewed as one-object enriched categories. Let $f:\mathcal{P}\to \thin{\mathcal{P}}$ be the unique functor where $\thin{\mathcal{P}}$ is the enriched small category which has one object and the identity map as the only map. Then the functor \[f^*:[\thin{\mathcal{P}},\top_Q]^{proj}_0\to [\mathcal{P},\top_Q]^{proj}_0\] is a right Quillen equivalence.
\eth

\bpf The first step is to replace the projective model structures by the projective mixed model structures. Indeed, it suffices by Theorem~\ref{proj-mixed-model-structure-enriched-diagram} to prove that the Quillen adjunction \[f_*\dashv f^*:[\thin{\mathcal{P}},\top_m]^{proj}_0 \leftrightarrows [\mathcal{P},\top_m]^{proj}_0\] is a Quillen equivalence to complete the proof. It then suffices to prove that the total derived functors induce an equivalence of categories between the homotopy categories by \cite[Proposition~1.3.13]{MR99h:55031}. To prove this fact, we can work in the more general setting of enriched homotopical categories in the sense of \cite{shuhomotopycolim}. We want to apply \cite[Proposition~22.5]{shuhomotopycolim}. In the language of \cite{shuhomotopycolim}, we have to prove that the two enriched small categories $\thin{\mathcal{P}}$ and $\mathcal{P}$ are very good for the tensor structure generated by the binary product of $\top$. To prove the latter fact, we have to use \cite[Theorem~23.12]{shuhomotopycolim}. It is easy to check that the monoidal model category $(\top_m,\p)$ is simplicial and that the binary product gives rise to a Quillen two-variable enriched adjunctions. All spaces of maps of $\thin{\mathcal{P}}$ and $\mathcal{P}$ are cofibrant in $\top_m$. To check the hypotheses of \cite[Theorem~23.12]{shuhomotopycolim}, it remains to check that any injection of a singleton in one of the space of maps of $\mathcal{G}$ or $\mathcal{M}$ is actually a mixed cofibration (cf. \cite[Definition~23.11]{shuhomotopycolim}). It is precisely what is proved in Proposition~\ref{well-pointed}. 
\epf


\begin{thebibliography}{DHKS04}
	
	\bibitem[AR94]{TheBook}
	J.~Ad{\'a}mek and J.~Rosick{\'y}.
	\newblock {\em Locally presentable and accessible categories}.
	\newblock Cambridge University Press, Cambridge, 1994.
	
	\bibitem[Bor94a]{Borceux1}
	F.~Borceux.
	\newblock {\em Handbook of categorical algebra. 1}.
	\newblock Cambridge University Press, Cambridge, 1994.
	\newblock Basic category theory.
	
	\bibitem[Bor94b]{Borceux2}
	F.~Borceux.
	\newblock {\em Handbook of categorical algebra. 2}.
	\newblock Cambridge University Press, Cambridge, 1994.
	\newblock Categories and structures.
	
	\bibitem[BQR98]{EnrichedSketch}
	F.~Borceux, C.~Quinteiro, and J.~Rosick\'y.
	\newblock A theory of enriched sketches.
	\newblock {\em Theory Appl. Categ.}, 4:No. 3, 47--72, 1998.
	
	\bibitem[BR13]{Barthel-Riel}
	T.~Barthel and E.~Riehl.
	\newblock On the construction of functorial factorizations for model
	categories.
	\newblock {\em Algebr. Geom. Topol.}, 13(2):1089--1124, 2013.
	
	\bibitem[Cam19]{TimHelp}
	T.~Campion.
	\newblock The binary product of two presentable objects.
	\newblock MathOverflow, 2019.
	\newblock URL:https://mathoverflow.net/q/306129 (version: 2019-03-12).
	
	\bibitem[Col06]{mixed-cole}
	M.~Cole.
	\newblock Mixing model structures.
	\newblock {\em Topology Appl.}, 153(7):1016--1032, 2006.
	
	\bibitem[CS02]{monographie_hocolim}
	W.~Chach{\'o}lski and J.~Scherer.
	\newblock Homotopy theory of diagrams.
	\newblock {\em Mem. Amer. Math. Soc.}, 155(736):x+90, 2002.
	
	\bibitem[DHKS04]{HomotopicalCategory}
	W.~G. Dwyer, P.~S. Hirschhorn, D.~M. Kan, and J.~H. Smith.
	\newblock {\em Homotopy limit functors on model categories and homotopical
		categories}, volume 113 of {\em Mathematical Surveys and Monographs}.
	\newblock American Mathematical Society, Providence, RI, 2004.
	
	\bibitem[DI04]{hocolimfacile}
	D.~Dugger and D.~C. Isaksen.
	\newblock Topological hypercovers and {$\mathbb{A}^1$}-realizations.
	\newblock {\em Math. Z.}, 246(4):667--689, 2004.
	
	\bibitem[FR08]{FR}
	L.~Fajstrup and J.~Rosick\'y.
	\newblock A convenient category for directed homotopy.
	\newblock {\em Theory and Applications of Categories}, 21(1):pp 7--20, 2008.
	
	\bibitem[Gau03]{model3}
	P.~Gaucher.
	\newblock A model category for the homotopy theory of concurrency.
	\newblock {\em Homology, Homotopy and Applications}, 5(1):p.549--599, 2003.
	
	\bibitem[Gau05]{model2}
	P.~Gaucher.
	\newblock Comparing globular complex and flow.
	\newblock {\em New York Journal of Mathematics}, 11:p.97--150, 2005.
	
	\bibitem[Gau09]{mdtop}
	P.~Gaucher.
	\newblock Homotopical interpretation of globular complex by multipointed
	d-space.
	\newblock {\em Theory Appl. Categ.}, 22:No.\ 22, 588--621 (electronic), 2009.
	
	\bibitem[Gau18]{leftdetflow}
	P.~Gaucher.
	\newblock Flows revisited: the model category structure and its left
	determinedness.
	\newblock arXiv:1806.08197, 2018.
	
	\bibitem[GKR18]{GKR18}
	R.~Garner, M.~Kedziorek, and E.~Riehl.
	\newblock Lifting accessible model structures.
	\newblock arXiv:1802.09889, 2018.
	
	\bibitem[Gra03]{mg}
	M.~Grandis.
	\newblock Directed homotopy theory. {I}.
	\newblock {\em Cah. Topol. G\'eom. Diff\'er. Cat\'eg.}, 44(4):281--316, 2003.
	
	\bibitem[Hat02]{MR1867354}
	A.~Hatcher.
	\newblock {\em Algebraic topology}.
	\newblock Cambridge University Press, Cambridge, 2002.
	
	\bibitem[Hir03]{ref_model2}
	P.~S. Hirschhorn.
	\newblock {\em Model categories and their localizations}, volume~99 of {\em
		Mathematical Surveys and Monographs}.
	\newblock American Mathematical Society, Providence, RI, 2003.
	
	\bibitem[Hov99]{MR99h:55031}
	M.~Hovey.
	\newblock {\em Model categories}.
	\newblock American Mathematical Society, Providence, RI, 1999.
	
	\bibitem[Kar18]{2948666}
	A.~Karagila.
	\newblock Regular cardinal for a fixpoint of the exponential.
	\newblock Mathematics Stack Exchange, 2018.
	\newblock URL:https://math.stackexchange.com/q/2948666 (version: 2018-10-09).
	
	\bibitem[Kel82]{Kelleyb}
	G.~M. Kelly.
	\newblock Structures defined by finite limits in the enriched context. {I}.
	\newblock {\em Cahiers Topologie G\'{e}om. Diff\'{e}rentielle}, 23(1):3--42,
	1982.
	\newblock Third Colloquium on Categories, Part VI (Amiens, 1980).
	
	\bibitem[Kel05]{KellyEnriched}
	G.~M. Kelly.
	\newblock Basic concepts of enriched category theory.
	\newblock {\em Repr. Theory Appl. Categ.}, (10):vi+137 pp. (electronic), 2005.
	\newblock Reprint of the 1982 original [Cambridge Univ. Press, Cambridge;
	MR0651714].
	
	\bibitem[Law18]{311786}
	T.~Lawson.
	\newblock Homotopy type of a specific discrete monoid.
	\newblock MathOverflow, 2018.
	\newblock URL:https://mathoverflow.net/q/311786 (version: 2018-10-01).
	
	\bibitem[McD79]{McDuff}
	D.~McDuff.
	\newblock On the classifying spaces of discrete monoids.
	\newblock {\em Topology}, 18(4):313--320, 1979.
	
	\bibitem[ML98]{MR1712872}
	S.~Mac~Lane.
	\newblock {\em Categories for the working mathematician}.
	\newblock Springer-Verlag, New York, second edition, 1998.
	
	\bibitem[Mos19]{MoserLyne}
	L.~Moser.
	\newblock Injective and projective model structures on enriched diagram
	categories.
	\newblock {\em Homology, Homotopy and Applications}, 21(2):pp.279--300, 2019.
	
	\bibitem[Pia91]{PiaDgr}
	R.~J. Piacenza.
	\newblock Homotopy theory of diagrams and {CW}-complexes over a category.
	\newblock {\em Canad. J. Math.}, 43(4):814--824, 1991.
	
	\bibitem[Rap10]{Raptis-Strom}
	G.~Raptis.
	\newblock Homotopy theory of posets.
	\newblock {\em Homology Homotopy Appl.}, 12(2):211--230, 2010.
	
	\bibitem[Rez02]{REZK200265}
	C.~Rezk.
	\newblock Every homotopy theory of simplicial algebras admits a proper model.
	\newblock {\em Topology and its Applications}, 119(1):65 -- 94, 2002.
	
	\bibitem[Ros09]{MR2506258}
	J.~Rosick{\'y}.
	\newblock On combinatorial model categories.
	\newblock {\em Appl. Categ. Structures}, 17(3):303--316, 2009.
	
	\bibitem[Ros17]{MR3638359}
	J.~Rosick\'{y}.
	\newblock Accessible model categories.
	\newblock {\em Appl. Categ. Structures}, 25(2):187--196, 2017.
	
	\bibitem[Shu09]{shuhomotopycolim}
	M.~Shulman.
	\newblock Homotopy limits and colimits and enriched homotopy theory.
	\newblock arXiv:math/0610194, 2009.
	
	\bibitem[Str66]{MR0211403}
	A.~Str{\o}m.
	\newblock Note on cofibrations.
	\newblock {\em Math. Scand.}, 19:11--14, 1966.
	
	\bibitem[Vog71]{MR45:9323}
	R.~M. Vogt.
	\newblock Convenient categories of topological spaces for homotopy theory.
	\newblock {\em Arch. Math. (Basel)}, 22:545--555, 1971.
	
\end{thebibliography}
\end{document}